\documentclass[twoside,12pt,leqno]{article}
\advance\topmargin by -\headheight\advance\topmargin by -\headsep
 \addtocounter{page}{0}
\newcommand{\tmvolxx}{17}
\newcommand{\tmyearyyyy}{2007}

\newcommand{\FirstPageHead}[3]{
\vskip -8mm
\centerline {Travaux math\'ematiques, \quad
Volume #1 (#2),
#3,\quad \copyright\  Universit\'e du Luxembourg}}\vspace{-3mm}

\usepackage{amsmath}
\usepackage{amsthm}
\usepackage{amssymb}
\usepackage{amscd}
\usepackage{graphicx}
\usepackage{epsfig}
\input xy
\xyoption{all}

\numberwithin{equation}{section}
\newtheorem{theorem}{Theorem}[section]
\newtheorem{lemma}[theorem]{Lemma}
\newtheorem{proposition}[theorem]{Proposition}

\theoremstyle{definition}
\newtheorem{definition}[theorem]{Definition}
\newtheorem{example}[theorem] {Example}

\numberwithin{equation}{section}

\allowdisplaybreaks
\begin{document}
\thispagestyle{empty}
\FirstPageHead{\tmvolxx}{\tmyearyyyy}{\pageref{firstpage}--\pageref{lastpage}}
\label{firstpage}


\hyphenation{grou-po-ids}
\hyphenation{co-i-so-tro-pic}
\hyphenation{co-i-so-tro-pe}
\hyphenation{con-stra-i-nt}
\hyphenation{sub-ma-ni-fold}
\hyphenation{Mac-ken-zie}
\hyphenation{group-oid}
\numberwithin{equation}{section}

\newcommand{\hs}[1]{\hspace{#1cm}}
\newcommand{\vs}[1]{\vspace{#1cm}}

\newcommand{\pf}{\\\noindent{\footnotesize{PROOF.} }}
\newcommand{\pff}{\noindent{\footnotesize{PROOF.} }}
\newcommand{\prf}[1]{\spa{\footnotesize{PROOF OF {#1}}.} }

\newcommand{\fin}[1]{\hspace{-#1pt}\hs{30}$\square$}
\newcommand{\find}[1]{\hspace{-#1pt}\hs{30}$FD$\smallskip}
\def\spa{\vskip0.25cm\noindent}

\def\Dots{,\,\dots\,,}
\def\eps{\varepsilon}
\let\Hat=\widehat

\def\id{{\sf{id}}}
\def\ker{{\sf{ker}}}
\def\im{{\sf{Im}}}
\def\rank{{\sf{Rank}}}
\def\comp{\hs{-0.035}\circ\hs{-0.035}}
\newcommand{\hbxs}[1]{\hbox{\small{$#1$}}}
\def\fii{\varphi}
\def\inc{\hookrightarrow}
\def\sur{\rightrightarrow}


\def\pa{\partial}
\def\dd{\mathrm{d}}
\newcommand{\pad}[2]{\frac{\pa{#1}}{\pa{#2}}}
\newcommand{\tod}[2]{\frac{\dd{#1}}{\dd{#2}}}
\newcommand{\vad}[2]{\frac{\delta{#1}}{\delta{#2} }}
\newcommand{\ddt}[1]{\frac{\dd{#1}}{\dd{t}}}
\newcommand{\ddtpa}[1]{\frac{\pa{#1}}{\pa{x}}}
\def\cif{\calC^\infty}
\def\frax{\mathfrak{X}}
\def\ham{{\sf{Ham}}}

\def\osharp{\omega^\sharp}
\def\oosharp{\Omega^\sharp}
\def\pish{\pi^\sharp}
\def\ppish{\Pi^\sharp}
\def\ohm{\omega}
\def\Ohm{\Omega}

\def\frag{\mathfrak{g}}
\def\frah{\mathfrak{h}}
\def\frack{\mathfrak{k}}
\def\frai{\mathfrak{i}}
\def\frak{\mathfrak{c}}

\def\gpd{\,\lower1pt\hbox{$\rightarrow$}\hskip-4.2mm\raise3.5pt
               \hbox{$\rightarrow$}\,}
\def\dpg{\,\lower1pt\hbox{$\leftarrow$}\hskip-3.75mm\raise3pt
               \hbox{$\leftarrow$}\,}

\newcommand{\ra}{\rangle}
\newcommand{\la}{\langle}
\newcommand{\backl}{\mathbin{\vrule width1.5ex height.4pt\vrule height1.5ex}}
\newcommand{\per}{\backl}

\newcommand{\smalcirc}{\mbox{\tiny{$\circ $}}}

\newcommand{\poidd }[2]{#1\gpd #2}
\newcommand{\poiddd }[3]{ (#1\gpd #2, \alpha_{#3}, \beta_{#3})}

\newcommand{\Palg}{P(\tTsM)}
\newcommand{\Galg}{G(\tTsM)}
\newcommand{\PGamma}{P_0\Gamma(\tTsM)}

\newcommand{\brak}[2]{[{#1,#2}]}
\newcommand{\poib}[2]{\{{#1,#2}\}}

\newcommand{\CC}{\mathbb{C}}
\newcommand{\RR}{\mathbb{R}}
\newcommand{\QQ}{\mathbb{Q}}
\newcommand{\ZZ}{\mathbb{Z}}
\newcommand{\NN}{\mathbb{N}}

\newcommand{\HH}{\mathbb{H}}

\newcommand{\FF}{\mathbb{F}}

\newcommand{\Z}[1]{\mathbb{Z}_{#1}}

\newcommand{\ds}{\mathbb{S}}
\newcommand{\BB}{\mathbb{B}}
\newcommand{\DD}{\mathbb{D}}

\newcommand{\calA}{\mathcal{A}}
\newcommand{\calB}{\mathcal{B}}
\newcommand{\calC}{\mathcal{C}}
\newcommand{\calD}{\mathcal{D}}
\newcommand{\calF}{\mathcal{F}}
\newcommand{\calH}{\mathcal{H}}
\newcommand{\calV}{\mathcal{V}}
\newcommand{\calS}{\mathcal{S}}
\newcommand{\calM}{\mathcal{M}}
\newcommand{\calK}{\mathcal{K}}
\newcommand{\calG}{\mathcal{G}}
\newcommand{\calI}{\mathcal{I}}
\newcommand{\calL}{\mathcal{L}}
\newcommand{\call}{\ell}
\newcommand{\calO}{\mathcal{O}}
\newcommand{\calR}{\mathcal{R}}
\newcommand{\tcalG}{\Tilde{\mathcal{G}}}
\newcommand{\ttcalG}{\Tilde{\tcalG}}
\newcommand{\ucalG}{\underline{\mathcal{G}}}
\newcommand{\bucalG}{\boldsymbol{\ucalG}}
\newcommand{\calE}{\mathcal{E}}
\newcommand{\calP}{\mathcal{P}}

\let\sf=\mathsf
\let\tsf=\textsf
\newcommand{\etsf}[1]{\tsf{\emph{#1}}}

\newcommand{\sfA}{\mathsf{A}}
\newcommand{\sfB}{\mathsf{B}}
\newcommand{\sfH}{\mathsf{H}}
\newcommand{\sfS}{\mathsf{S}}
\newcommand{\sfC}{\mathsf{C}}
\newcommand{\sfM}{\mathsf{M}}
\newcommand{\sfK}{\mathsf{K}}
\newcommand{\sfG}{\mathsf{G}}
\newcommand{\sfI}{\mathsf{I}}
\newcommand{\sfL}{\mathsf{L}}
\newcommand{\sfO}{\mathsf{O}}
\newcommand{\tsfG}{\Tilde{\mathsf{G}}}
\newcommand{\ttsfG}{\Tilde{\tsfG}}
\newcommand{\usfG}{\underline{\mathsf{G}}}
\newcommand{\busfG}{\boldsymbol{\usfG}}
\newcommand{\sfE}{\mathsf{E}}
\newcommand{\sfP}{\mathsf{P}}

\let\rm=\mathrm
\let\trm=\textrm
\newcommand{\rmA}{\mathrm{A}}
\newcommand{\rmB}{\mathrm{B}}
\newcommand{\rmH}{\mathrm{H}}
\newcommand{\rmS}{\mathrm{S}}
\newcommand{\rmC}{\mathrm{C}}
\newcommand{\rmM}{\mathrm{M}}
\newcommand{\rmK}{\mathrm{K}}
\newcommand{\rmG}{\mathrm{G}}
\newcommand{\rmI}{\mathrm{I}}
\newcommand{\rmL}{\mathrm{L}}
\newcommand{\rmO}{\mathrm{O}}
\newcommand{\trmG}{\Tilde{\mathrm{G}}}
\newcommand{\ttrmG}{\Tilde{\trmG}}
\newcommand{\urmG}{\underline{\mathrm{G}}}
\newcommand{\burmG}{\boldsymbol{\urmG}}
\newcommand{\rmE}{\mathrm{E}}
\newcommand{\rmP}{\mathrm{P}}

\let\Bar=\overline
\let\ol=\overline
\let\ul=\underline

\let\bgn=\begin

\newcommand{\temise}[1]{{
\begin{itemize}
{#1}
\end{itemize}
}}
\def\it{\item[i)]}
\def\iit{\item[ii)]}
\def\iiit{\item[iii)]}
\def\ivt{\item[iv)]}
\def\vt{\item[v)]}
\def\vit{\item[vi)]}
\def\viit{\item[vii)]}
\def\viiit{\item[viii)]}
\def\punto{\item[$\cdot$]}

\newcommand{\be}{\begin{eqnarray*}}
\newcommand{\ee}{\end{eqnarray*}}

\newcommand{\Exasss}[4]{{\Exa{{#1}\temise{
\item[1.] {#2}
\item[3.] {#3}
\item[4.] {#4}
}}}}

\def\daga{^{\pmb{++}}}

\def\ttar{\mathrm{\hbox{\small{$\mathrm{t}$}}}}
\def\tar{\mathrm{t}}
\def\ssor{\mathrm{\hbox{\small{$\mathrm{s}$}}}}
\def\sor{\mathrm{s}}
\def\vup{\hbox{\tiny{${}_{V}$}}}
\def\hup{\hbox{\tiny{${}_{H}$}}}
\def\vuppp{\hbox{\tiny{${}^{3V}$}}}

\def\ttv{\,{\ttar}\vup}
\def\tsv{\,{\ttar}_{v}}
\def\stv{\,{\ssor}\vup}
\def\ssv{\,{\ssor}_{v}}

\def\sc{{\sor}^c}
\def\tc{{\tar}^c}

\def\tth{\,{\ttar}\hup}
\def\tsh{\,{\ttar}_{h}}
\def\sth{\,{\ssor}\hup}
\def\ssh{\,{\ssor}_{h}}

\def\nttv{{\mathrm{t}}\vup}
\def\ntsv{{\mathrm{t}}^{v}}
\def\nstv{{\mathrm{s}}\vup}
\def\nssv{{\mathrm{s}}^{v}}

\def\ntth{{\mathrm{t}}\hup}
\def\ntsh{{\mathrm{t}}^{h}}
\def\nsth{{\mathrm{s}}\hup}
\def\nssh{{\mathrm{s}}^{h}}

\def\etv{{\eps}\vup}
\def\esv{{\eps}^{v}}
\def\etv{{\eps}\vup}
\def\esv{{\eps}^{v}}
\def\ec{{\eps}^{c}}

\def\eth{{\eps}\hup}
\def\esh{{\eps}^{h}}
\def\eth{{\eps}\hup}
\def\esh{{\eps}^{h}}

\def\itv{{\iota}\vup}
\def\isv{{\iota}^{v}}
\def\itv{{\iota}\vup}
\def\isv{{\iota}^{v}}

\def\ith{{\iota}\hup}
\def\ish{{\iota}^{h}}
\def\ith{{\iota}\hup}
\def\ish{{\iota}^{h}}

\def\mtv{{\mu}\vup}
\def\msv{{\mu}^{v}}
\def\mtv{{\mu}\vup}
\def\msv{{\mu}^{v}}

\def\mth{{\mu}\hup}
\def\msh{{\mu}^{h}}
\def\mth{{\mu}\hup}
\def\msh{{\mu}^{h}}

\def\dtv{{\delta}\vup}
\def\dsv{{\delta}^{v}}
\def\dtv{{\delta}\vup}
\def\dsv{{\delta}^{v}}

\def\dth{{\delta}\hup}
\def\dsh{{\delta}^{h}}
\def\dth{{\delta}\hup}
\def\dsh{{\delta}^{h}}

\def\epsh{\hat{\eps}}
\def\sh{\hat{\sor}}
\def\th{\hat{\tar}}
\def\iotah{\hat{\iota}}
\def\muh{\hat{\mu}}

\def\epst{\tilde{\eps}}
\def\st{\tilde{s}}
\def\tt{\tilde{t}}
\def\iotat{\tilde{\iota}}
\def\mut{\tilde{\mu}}

\newcommand{\fib}[2]{\,{}_{#1}\hs{-0.05}\times\hs{-0.05}{}_{#2}\,}
\newcommand{\Lie}[1]{\mathsf{Lie}[#1]}
\newcommand{\lie}[2]{\mathsf{Lie}_{#1}[#2]}

\def\lagpd{\mathcal{LA}\hbox{-groupoid}}
\def\pgpd{(\calP,\Pi)\gpd M}
\def\pgpdb{(\ol{\calP},{\ol\Pi})\gpd M}
\def\poidm{\calP\gpd M}
\def\poidmb{\bar{\calP}\gpd M}
\def\bagd{(A,A^*)\rightarrow M}

\def\tsp{T^*\calP}
\newcommand{\gr}[1]{\sf{\Gamma}({\!#1})}
\newcommand{\ki}[2]{\mathfrak{X}^{#1}{(#2)}}




\markboth{Luca Stefanini}{On the integration of
$\mathcal{LA}$-groupoids and duality for Poisson groupoids}
$ $
\bigskip

\bigskip

\centerline{{\Large On the integration of
$\mathcal{LA}$-groupoids and duality for Poisson groupoids}}

\bigskip
\bigskip
\centerline{{\large by
Luca Stefanini\footnote{The author acknowledges support by SNF-grant Nr.20-113439.}}}

\vspace*{.7cm}

\begin{abstract}
In this note a functorial approach to the integration problem of an
$\mathcal{LA}$-groupoid to a double Lie groupoid is discussed. To do that, we study the notions of
fibered products in the categories of Lie groupoids and Lie algebroids, giving
necessary and sufficient conditions for the existence of such. In particular, it
turns out, that the fibered product of Lie algebroids along integrable morphisms
is always integrable by a fibered product of Lie groupoids. We show that to every
$\lagpd$ with
integrable top structure one can associate a differentiable graph in the category
of Lie groupoids, which is an integrating double Lie groupoid, whenever some
lifting conditions for suitable Lie algebroid homotopies are fulfilled; the result specializes to the case of a Poisson
groupoid, yielding a symplectic double groupoid, provided our conditions on the
associated $\lagpd$ are satisfied\footnote{\emph{2000 Mathematics Subject Classification}:
Primary 58H05; Secondary 17B66, 18D05, 22A22.}.
\end{abstract}
\pagestyle{myheadings}
\section*{Introduction}
In recent years, fundamental questions in Lie theory for Lie algebroids and Lie
groupoids have been answered; namely, optimal generalizations of Lie's theorems
have been discovered. Examples of non-integrable Lie algebroids
already appeared in \cite{AM} and the problem to find general integrability
conditions has been standing for a long time.\\
The $-$quite non-trivial$-$ theory of morphisms of Lie algebroids was developed by
Higgins and Mackenzie in \cite{HM}. Later on Mackenzie and Xu proved \cite{MX2} that
morphisms of integrable Lie algebroids are integrable to morphisms of Lie
groupoids, provided the domain groupoid has 1-connected source fibres. An
independent proof by Moerdijk and Mr\v cun appeared in \cite{MM}. In the same paper
 the authors also show
 that to every source connected Lie groupoid one can
associate a unique source 1-connected ``cover'' with the same Lie algebroid; moreover they
prove that every Lie subalgebroid of an integrable Lie algebroid $A$ is integrable by an
(only) immersed subgroupoid of the source 1-connected integration of $A$. The source
1-connected cover of a source connected Lie groupoid $\calG$ is obtained as the
quotient of the monodromy groupoid $\sf{Mon}(\calG,\sor)$ associated with
the source foliation on $\calG$ with respect to the natural action by
right translation of $\calG$ itself.
 It turns out \cite{CrFr}, that the Lie groupoid
$\sf{Mon}(\calG,\sor)/\calG$ can be equivalently described as the quotient of the so called
$\calG$-paths, paths along the source fibers starting from the base manifold, by
homotopy within the source fibers, relative to the end points. Crainic and
Fernandes showed that both the notions of $\calG$-paths and their homotopy can be
characterized in terms of the Lie algebroid $A$ of $\calG$; namely,
$\calG$-paths are in bijective correspondence with $A$-paths,
i.e. morphisms of Lie algebroids $TI\rightarrow A$, and $\calG$-paths are
homotopic iff the corresponding $A$-paths
are $A$-homotopic, being $A$-homotopopies morphisms of Lie algebroids $TI^{\times 2}\rightarrow A$,
satisfying suitable boundary conditions. The quotient
$\mathcal{W}(A):=\{A\hbox{-}paths\}/A\hbox{-}homotopy$,
a.k.a. the \emph{Weinstein groupoid},
carries a natural groupoid structure, induced, roughly,
by concatenation of paths; it is always a topological groupoid and Crainic and
Fernandes finally delivered a necessary and sufficient integrability condition for
Lie algebroids, which is to be understood as the obstruction to put a smooth
structure on the associated Weinstein groupoid.\\
The construction of the Weinstein groupoid was anticipated by Cattaneo and Felder
\cite{CtFl}
in the special case of the Lie algebroid of a Poisson manifold. Their approach
involves the symplectic reduction of the phase space of the Poisson sigma model
and yields the symplectic groupoid of the target Poisson manifold,
in the integrable case.

In this paper we study (part of) the categorified version of this story.
Ehresmann's categorification of a groupoid is a groupoid(-object) in the category of
groupoids; this is a symmetric notion and it makes sense to regard such a structure as a
``double groupoid''.  A double Lie groupoid is, essentially, a ``Lie groupoid in the category
of Lie groupoids''; one can apply the Lie functor to the object of a double Lie
groupoid, to obtain an $\lagpd$, i.e. a  ``Lie groupoid in the category of
Lie algebroids''. The application of the Lie functor can
still be iterated; the result, a double Lie algebroid, is the best approximation to
what one would mean as a ``Lie algebroid in the category of Lie algebroids''. Such
double structures do arise in nature, especially from Poisson geometry and
the theory of
Poisson actions.

After reviewing the main definitions and known integrability results related to Lie
bialgebroids and Poisson groupoids (\S\ref{lapo}), we address the integration
problem of an $\lagpd$ to a double Lie groupoid. In \S\ref{3} we study the fibred
products of Lie algebroids and Lie groupoids. We show that, whenever a fibered product
of Lie algebroids exists as a vector bundle, it carries a unique natural Lie algebroid
structure; the analogous property does not hold for Lie groupoids. We find,
however a
necessary and sufficient transversality condition for fibered products of Lie groupoids
to stay in the category. We develop our integration approach in \S\ref{fungo},
also considering the case of the $\lagpd$ of a Poisson groupoid in relation with
duality issues (\S\ref{pioppo}). We introduce differentiable graphs with structure, Lie groupoids
``without a multiplication'', and show
that an $\lagpd$ with integrable top Lie algebroid is always integrable to an invertible
unital graph in the category of Lie groupoids; moreover, we prove a natural integrability
result for fibered products of integrable Lie algebroids. As a consequence the
integrating graph of an $\lagpd$ can be endowed with a further compatible multiplication
making it a double Lie groupoid, under some connectivity assumptions on its second and
third nerve. Surprisingly, the connectivity assumptions are implied by some
suitable lifting conditions for Lie algebroid paths and Lie algebroid homotopies, depending only
on the original $\lagpd$. We shall remark, however,
that our requirements are far from being
necessary integrability conditions and appear
quite restrictive.\\
Lastly, we comment on an alternative approach to the $-$2 steps in 1$-$
integration of a Lie bialgebroid to a symplectic double groupoid, within the
framework of symplectic reduction of the Courant sigma model.
\subsection*{Notations and conventions} We 
denote with $\sor$ and $\tar$
the source and target maps
of a Lie groupoid, with $\eps$ the unit section, with $\iota$ the inversion and
with $\mu$ the partial
multiplication. The anchor of a Lie algebroid is typically denoted with $\rho$.
Nowhere in this paper
$\ol{P}$ denotes the opposite Poisson structure on a Poisson manifold $P$.
``Fibered product'' is meant as a categorical pullback; with pullback it is
meant ``pullback along a map'', such as the vector bundle pullback.
We shall denote with $\sf\Gamma(f)$ the graph of a map $f:M\to N$ and
regard it as a fibered product
$M\fib{f}{}N\equiv M\fib{f}{\id_N}N$. If two smooth maps
$f_{1,2}:M^{1,2}\to N$ are transversal, we shall write $f_1\pitchfork f_2$. For
two vector bundle maps $\phi_{1,2}:E^{1,2}\to F$ to be transversal,
$\phi_1\pmb{\pitchfork}\phi_2$, means that $\phi_1$ and $\phi_2$ are transversal as
smooth maps, so are the corresponding base maps $f_{1,2}:M_{1,2}\to N$ and
$
\phi_1-\phi_2:\left.(E^1\times
E^2)\right|_{M^1\fib{f_1}{f_2}M^2}\rightarrow F
$
has constant rank, so that the
fibered product $E^1\fib{\phi_1}{\phi_2}E^2$ carries a vector bundle structure over
the fibered product $M^1\fib{f_1}{f_2}M^2$.
%
%
%
%
\subsection*{Acknowledgments}
First of all, I wish to thank Alberto Cattaneo,
for helping me glimpse how deep the rabbit
hole goes. I am very grateful to Kirill Mackenzie for his advice and
encouragement. Special thanks go to
Ping Xu, for many stimulating conversations and
the warmest hospitality down at State College.\\
Moreover, I wish to thank the anonymous referee, for valuable comments
and for suggesting interesting possible directions of further research,
and  Kirill Mackenzie, also
for carefully reading this paper, while it was being refereed, and
for providing useful suggestions to improve the presentation.%
\section{$\lagpd$s, (symplectic) double Lie groupoids and Poisson
groupoids}\label{lapo}
A double Lie groupoid is a Lie groupoid object in the
category of Lie groupoids.
\begin{definition}[\cite{M1}]\label{dlg}
A \tsf{\textsf{double Lie groupoid}}
$\sf{D}:=(\calD,\calH,\calV,M)$
$$
\xy
*+{}="0",    <-0.7cm,0.7cm>
*+{\calD}="1", <0.7cm,0.7cm>
*+{\calV}="2", <-0.7cm,-0.7cm>
*+{\calH}="3", <0.7cm,-0.7cm>
*+{M}="4",
\ar  @ <0.07cm>   @{->} "1";"2"^{}
\ar  @ <-0.07cm>  @{->} "1";"2"_{}
\ar  @ <-0.07cm>  @{->} "1";"3"_{}
\ar  @ <0.07cm>   @{->} "1";"3"^{}
\ar  @ <0.07cm>   @{->} "2";"4"^{}
\ar  @ <-0.07cm>  @{->} "2";"4"_{}
\ar  @ <-0.07cm>  @{->} "3";"4"_{}
\ar  @ <0.07cm>   @{->} "3";"4"^{}
\endxy
$$
is a groupoid object in the category of groupoids (i.e. a \tsf{double
groupoid} in the sense of Ehresmann), such that
$\poidd{\calD}{\calV}$, $\poidd{\calD}{\calH}$, $\poidd{\calH}{M}$,
$\poidd{\calV}{M}$ are Lie groupoids
and
the double
source map
$$\mathbb{S}\doteq(\nsth,\nstv):\calD\rightarrow \calH\fib{\ssh}{\ssv}\calV$$
is
submersive\footnote{In \cite{M1} the double source map is required to be also
surjective; this condition does not really play a r\^ole in the study of the
internal structure of a double Lie groupoid and the descent to double
Lie algebroids. Moreover, there are interesting examples, such as Lu and
Weinstein's double of a Poisson group (\ref{lwd}) for instance, which do not
fulfill the double source surjectivity condition.}.
\end{definition}
The definition is symmetric and the total space of a double Lie groupoid
can be regarded as a
groupoid object either horizontally or vertically; the groupoid vertical, resp.
horizontal, structural maps (unit section, source, target, inversion and multiplication)
are morphisms of Lie groupoids for the vertical, resp. horizontal, structures.
Note that the submersivity condition on the double source map makes the
domains of the top multiplications Lie groupoids (see proposition (\ref{Fibbia}) for a
justification of this fact).

Applying the Lie functor  horizontally, or vertically, yields an
$\mathcal{LA}$-groupoid.
\begin{definition}[\cite{M1}]
An $\mathcal{LA}$-\tsf{groupoid} $\sf{\Ohm}:=(\Ohm,A,\calG,M)$
$$
\xy
*+{}="0",    <-0.7cm,0.7cm>
*+{\Omega}="1", <0.7cm,0.7cm>
*+{\calG}="3", <-0.7cm,-0.7cm>
*+{A}="2", <0.7cm,-0.7cm>
*+{M}="4",
\ar  @ <-0.07cm>   @{->} "1";"2"_{\Hat{}\quad}
\ar  @ <0.07cm>    @{->} "1";"2"^{}
\ar            @{->} "1";"3"^{}
\ar                @{->} "2";"4"_{}
\ar  @ <0.07cm>    @{->} "3";"4"^{}
\ar  @ <-0.07cm>   @{->} "3";"4"_{}
\endxy
$$
is a groupoid object in the category of Lie algebroids, such that
$\poidd{\Omega}{A}$ and $\poidd{\calG}{M}$ are
Lie groupoids and
the double source map
$$
\$\doteq(\sh,\rm{Pr}):
\Omega\rightarrow A\fib{\rm{pr}}{\sor}\calG$$
is a surjective submersion.
\end{definition}
An $\lagpd$ is a double groupoid for the horizontal additive
groupoids and the double source map should be understood with respect to this structure. As a direct
consequence from the definition, the vector bundle projection, zero section, fibrewise addition and
scalar multiplication of $\Ohm\rightarrow\calG$ are morphisms of Lie groupoids over the corresponding
maps of $A\rightarrow M$. Note that there is no natural way of characterizing
the Lie algebroid bracket of
$\Ohm$ as a morphism of Lie algebroids over the bracket of $A$. On the other hand,
an $\lagpd$ is indeed a Lie groupoid in the category of Lie algebroids.
\begin{example}\label{exa}
The typical examples are: for any Lie groupoid $\poidd{\calG}{M}$
$$\qquad\qquad
\begin{array}{ccc}
\xy
*+{}="0",    <-1.4cm,0.7cm>
*+{\calG\times\calG}="1", <0.4cm,0.7cm>
*+{\calG}="2", <-1.4cm,-0.7cm>
*+{M\times M}="3", <0.4cm,-0.7cm>
*+{M}="4",
\ar  @ <0.07cm>   @{->} "1";"2"^{}
\ar  @ <-0.07cm>  @{->} "1";"2"_{}
\ar  @ <-0.07cm>  @{->} "1";"3"_{}
\ar  @ <0.07cm>   @{->} "1";"3"^{}
\ar  @ <0.07cm>   @{->} "2";"4"^{}
\ar  @ <-0.07cm>  @{->} "2";"4"_{}
\ar  @ <-0.07cm>  @{->} "3";"4"_{}
\ar  @ <0.07cm>   @{->} "3";"4"^{}
\endxy
&\qquad\qquad&
\xy
*+{}="0",    <-0.7cm,0.7cm>
*+{T\calG}="1", <0.7cm,0.7cm>
*+{\calG}="3", <-0.7cm,-0.7cm>
*+{TM}="2", <0.7cm,-0.7cm>
*+{M}="4",
\ar  @ <-0.07cm>   @{->} "1";"2"_{}
\ar  @ <0.07cm>    @{->} "1";"2"^{}
\ar            @{->} "1";"3"^{}
\ar                @{->} "2";"4"_{}
\ar  @ <0.07cm>    @{->} "3";"4"^{}
\ar  @ <-0.07cm>   @{->} "3";"4"_{}
\endxy
\\
\end{array}\qquad\qquad,
$$
where the top vertical groupoid is a direct product and the horizontal groupoids
are pair groupoids, in the first case, while in the second, the top groupoid is
the tangent prolongation (each structural map is the tangent of the
corresponding map of $\calG$).
\end{example}
Applying the Lie functor once more yields a \tsf{double Lie algebroid} (See
\cite{M1,M7} for a definition). The notions of morphisms and sub-objects
of double structures are the
obvious ones.
\spa

Lie theory ``from double Lie groupoids to double Lie algebroids'' has been developed to a
satisfactory extent in recent years by Mackenzie \cite{M1}-\cite{M5}, \cite{M7}.
Integrability results for
double Lie algebroids to
$\lagpd$s and for $\lagpd$s to double Lie groupoids are known only for a restricted class of
examples arising from Poisson geometry, the main ones we are about to sketch.
\begin{definition}[\cite{W}]\label{pgrpd}
A \tsf{Poisson groupoid} is a Lie Groupoid $\poidm$ endowed with a
compatible
Poisson structure $\Pi\in\ki{2}{\calP}$, i.e. such that the graph
$\gr{\,\mu}\subset\calP^{\times 3}$ of the groupoid
multiplication $\mu$ is coisotropic with respect to the Poisson structure
$\Pi\times\Pi\times
-\Pi$.
\end{definition}
A Poisson groupoid with a non-degenerate Poisson structure is a
symplectic groupoid in the usual sense (i.e. $\gr{\,\mu}$ is Lagrangian, by
counting dimensions). For any Poisson groupoid $\pgpd$ \cite{W}:\smallskip\\
($i$) The unit section $\eps:M\inc\calP$ is a closed coisotropic embedding;
\smallskip\\
($ii$) The inversion map $\iota:\calP\rightarrow\calP$ is an anti-Poisson diffeomorphism;
\smallskip\\
($iii$) The source
invariant functions and the target
invariant functions define commuting
anti-isomorphic Poisson subalgebras of $\cif(\calP)$.\smallskip\\
As an easy consequence of property ($iii$) above, the base manifold of a Poisson
groupoid carries a unique Poisson structure making the source map
Poisson and the target map anti-Poisson. A symplectic groupoid provides a
\emph{symplectic realization}
of the Poisson structure induced
on the base (see, for example, \cite{CDW,CrFr2} for an account on symplectic
realizations and symplectic groupoids).

It turns out \cite{MX1}, that a Lie groupoid $\poidd{\calP}{M}$ with a Poisson structure $\Pi$
is a Poisson groupoid iff
\bgn{equation}\label{@}
\xy
*+{}="0",    <-0.7cm,0.7cm>
*+{\tsp}="1", <0.7cm,0.7cm>
*+{T\calP}="3", <-0.7cm,-0.7cm>
*+{A^*}="2", <0.7cm,-0.7cm>
*+{TM}="4",
\ar  @ <-0.07cm>   @{->} "1";"2"_{ }
\ar  @ <0.07cm>    @{->} "1";"2"^{ }
\ar            @{->} "1";"3"^{\Pi^\sharp}
\ar                @{->} "2";"4"_{}
\ar  @ <0.07cm>    @{->} "3";"4"^{}
\ar  @ <-0.07cm>   @{->} "3";"4"_{}
\endxy\qquad
\hbox{is a morphism of groupoids.}
\end{equation}
The base map in the above diagram is the restriction of the Poisson anchor to
$N^*M=\tsf{Ann}_{\tsp}TM$, which is to be canonically identified with the dual
bundle
to the Lie algebroid $A$ of $\calP$. The cotangent prolongation
groupoid \cite{CDW} $\poidd{T^*\calG}{A^*}$, can be defined for any Lie groupoid
$\poidd{\calG}{M}$
and it is the symplectic groupoid of the fibrewise linear Poisson structure
induced from the Lie algebroid $A$ of $\calG$ on $A^*$. If $\calP$ is a Poisson groupoid $A^*$, as the
conormal bundle to a coisotropic submanifold, carries a Lie algebroid structure
over $M$ and
\bgn{equation}\label{@@}
\xy
*+{}="0",    <-0.7cm,0.7cm>
*+{\tsp}="1", <0.7cm,0.7cm>
*+{\calP}="3", <-0.7cm,-0.7cm>
*+{A^*}="2", <0.7cm,-0.7cm>
*+{M}="4",
\ar  @ <-0.07cm>   @{->} "1";"2"_{ }
\ar  @ <0.07cm>    @{->} "1";"2"^{ }
\ar            @{->} "1";"3"^{}
\ar                @{->} "2";"4"_{}
\ar  @ <0.07cm>    @{->} "3";"4"^{}
\ar  @ <-0.07cm>   @{->} "3";"4"_{}
\endxy\qquad\qquad
\hbox{is an $\lagpd$.}\qquad\quad
\end{equation}
In fact (\ref{@}) is the compatibility condition between horizontal anchors and
vertical Lie groupoids. The compatibility with the Lie algebroid brackets can be
shown as a consequence of the duality between $\mathcal{PVB}$-groupoids and
$\lagpd$s (see \cite{M2} for details).

Recall that, if $\pgpd$ is a Poisson groupoid $(A,A^*)$, is a \tsf{Lie
bialgebroid} \cite{MX1}; that is, the Lie algebroid structures on $A$
and $A^*$ are compatible, in the sense that \cite{K} $(\Gamma(\wedge^\bullet A^*),
\,\wedge\,,\,\dd_A\,,\, [\:,\:]_{A^*})$
is a differential Gerstenhaber algebra for the Lie algebroid differential $\dd_A$
induced by $A$ and the graded Lie bracket $[\:,\:]_{A^*}$ on
$\Gamma(\wedge^\bullet A^*)[1]$ induced by $A^*$. The notion of Lie bialgebroid
is self dual ($(A,A^*)$ is a Lie bialgebroid iff so is $(A^*,A)$) and the
\tsf{flip} $(A^*,\bar{A})$ of a Lie bialgebroid (invert signs of the anchor and
bracket of $A$) is also a Lie bialgebroid. This leads to a notion of duality for
Poisson groupoids, essentially introduced in \cite{W}.
\begin{definition}[\cite{M2}]
Poisson groupoids $\poidd{(\calP_\pm,\Pi_\pm)}{M}$ are in \tsf{weak duality}
if the Lie bialgebroid of $\calP_+$ is isomorphic to the flip of the
Lie bialgebroid of $\calP_-$.
\end{definition}
There is an important integrability result for Lie bialgebroids.
\begin{theorem}[\cite{MX2}]\label{intb}
For any Lie bialgebroid $\bagd$, with $A$ integrable, there exists a unique
Poisson structure $\Pi\in\ki{2}{\calP}$ on the source 1-connected Lie groupoid $\poidm$ of $A$,
such that\smallskip\\
1.
$\pgpd$ is a Poisson groupoid,\smallskip\\
2.
The Lie algebroid on $A^*$ coincides with that
induced by $\Pi$.
\end{theorem}
Last result, in view of \cite{M2,M3,M7}, can be interpreted as an integrability
result for the cotangent double (Lie algebroid) of a Lie bialgebroid
$$
\xy
*+{}="0",    <-0.7cm,0.7cm>
*+{T^*A^*}="1", <0.7cm,0.7cm>
*+{A}="3", <-0.7cm,-0.7cm>
*+{A^*}="2", <0.7cm,-0.7cm>
*+{M}="4",
\ar     @{->} "1";"2"_{ }
\ar     @{->} "1";"3"^{}
\ar     @{->} "2";"4"_{}
\ar     @{->} "3";"4"_{}
\endxy
$$
to the $\lagpd$ (\ref{@@}) of the corresponding source 1-connected Poisson groupoid.
See also \cite{M6} for another example of an integrable double Lie algebroid arising from the Poisson action
of a Poisson group.\\
An easy consequence of theorem (\ref{intb}) is the following.
\begin{lemma}
Every integrable Poisson groupoid $\pgpd$ has a unique source 1-connected
weak dual Poisson groupoid $\poidd{(\ol{\calP}, \ol{\Pi})}{M}$.
\end{lemma}
\begin{proof}
Recall from \cite{MM} that a Lie subalgebroid of an integrable Lie algebroid
is integrable. Since $\calP$ is an integrable Poisson manifold $T^*\calP$ is an
integrable Lie algebroid and so is $A^*$ (embed $A^*$ in $T^*\calP$ using the
unit section of the cotangent prolongation groupoid). Apply theorem (\ref{intb})
 to
the flip of the Lie bialgebroid of $\calP$.
\end{proof}
A stronger notion of duality for Poisson groupoids arises from double groupoids.
\begin{definition}[\cite{M2}]
A \tsf{symplectic double groupoid} is a double Lie groupoid, whose total space
is endowed with a symplectic form, which is compatible with the top horizontal
and vertical Lie groupoids.
\end{definition}
It follows \cite{M2}, that the Poisson structures on the total spaces of the side
groupoids induced by the top horizontal and vertical groupoids are in weak duality.
This motivates the following definition.
\begin{definition}\hs{-0.2}\footnote{This notion was suggested to the author
by K. Mackenzie in a private discussion (2005).}
Two Poisson groupoids $\poidd{(\calP_\pm,\Pi_\pm)}{M}$
are in
\tsf{strong duality} if there exists a symplectic double groupoid
$$
\xy
*+{}="0",    <-0.7cm,0.7cm>
*+{\calS}="1", <0.7cm,0.7cm>
*+{\calP_+}="2", <-0.7cm,-0.7cm>
*+{\calP_-}="3", <0.7cm,-0.7cm>
*+{M}="4",
\ar  @ <0.07cm>   @{->} "1";"2"^{}
\ar  @ <-0.07cm>  @{->} "1";"2"_{}
\ar  @ <-0.07cm>  @{->} "1";"3"_{}
\ar  @ <0.07cm>   @{->} "1";"3"^{}
\ar  @ <0.07cm>   @{->} "2";"4"^{}
\ar  @ <-0.07cm>  @{->} "2";"4"_{}
\ar  @ <-0.07cm>  @{->} "3";"4"_{}
\ar  @ <0.07cm>   @{->} "3";"4"^{}
\endxy
$$
with the given groupoids as side structures, such that the symplectic form
induces the given Poisson structures. The \tsf{double of an integrable Poisson groupoid}
is a symplectic double groupoid $\calS$ realizing a strong
duality of $\calP$ and its weak dual Poisson groupoid $\ol{\calP}$.
\end{definition}
Thus, strongly dual Poisson groupoids admit a simultaneous integration (and
symplectic realization); the natural questions
to answer are:\smallskip\\
$\cdot$ \emph{Do Poisson groupoids integrate to symplectic double groupoids}?
\smallskip\\
And, more generally,
\smallskip\\
$\cdot$ \emph{Does weak duality imply strong duality for integrable Poisson
groupoids}?\smallskip\\
A positive answer to the first question was given in \cite{LW} by Lu and Weinstein
in the case of integrable Poisson groups, and \cite{LP} by Li and Parmentier in the case of a
class of coboundary dynamical Poisson groupoids.
\begin{example}\cite{LW}\label{lwd} Let $(\frag,\frag^*)$ be the tangent Lie bialgebra of a 1-connected Poisson group $G$
and
$\ol{G}$ be
the weak dual. The sum $\mathfrak{d}=\frag\oplus\frag^*$ carries a natural Lie algebra structure, obtained by
a double twist of the brackets on $\frag$ and $\frag^*$, the
\emph{Drinfel'd double} of $(\frag,\frag^*)$;
 let $D$ be the 1-connected
integration of such. Denote with $\lambda:G\hookrightarrow D$, resp.
$\rho:\ol{G}\hookrightarrow D $  the integrations of $\frag\hookrightarrow\mathfrak{d}$, resp.
$\frag^*\hookrightarrow\mathfrak{d}$. One can show that $D$ has a compatible
Poisson structure $\pi_D$, which happens to be non-degenerate on the submanifold of elements $d$, admitting a
decomposition $d=\lambda(g_+)\rho(\ol{g}_+)=\rho(\ol{g}_-)\lambda(g_-)$, $g_\pm\in G$ and
$\ol{g}_\pm\in\ol{G}$. Moreover there is also a natural double Lie groupoid
$$\qquad\qquad\qquad\qquad\qquad\qquad\qquad
\xy
*+{}="0",    <-0.7cm,0.7cm>
*+{\calD}="1", <0.7cm,0.7cm>
*+{G}="2", <-0.7cm,-0.7cm>
*+{\ol{G}}="3", <0.7cm,-0.7cm>
*+{\bullet}="4",
\ar  @ <0.07cm>   @{->} "1";"2"^{}
\ar  @ <-0.07cm>  @{->} "1";"2"_{}
\ar  @ <-0.07cm>  @{->} "1";"3"_{}
\ar  @ <0.07cm>   @{->} "1";"3"^{}
\ar  @ <0.07cm>   @{->} "2";"4"^{}
\ar  @ <-0.07cm>  @{->} "2";"4"_{}
\ar  @ <-0.07cm>  @{->} "3";"4"_{}
\ar  @ <0.07cm>   @{->} "3";"4"^{}
\endxy\qquad\qquad\calD=G\times\ol{G}\times\ol{G}\times G\quad.
$$
It turns out that the double subgroupoid, whose total space is
{\small$\calS=\{(g_+,\ol{g}_+,\ol{g}_-,g_-)\:|$\\ $\:\lambda(g_+)\rho(\ol{g}_+)
=\rho(\ol{g}_-)\lambda(g_-)\}$}
carries a compatible symplectic form, inducing the Poisson structures on $G$ and
$\ol{G}$, which is the inverse of the pullback of $\pi_D$, under
the natural local diffeomorphism $\calS\rightarrow D$.
Note that $\calS$ is, in general, neither vertically,
nor horizontally source (1-)connected.
\end{example}
\section{Fibered products in the categories of Lie algebroids and Lie
groupoids}\label{3}
Recall from \cite{HM} that, given Lie algebroids $A^{1,2}\rightarrow M^{1,2}$, with anchors $\rho_{1,2}$ and
brackets $[\cdot,\cdot]_{1,2}$ a \textsf{ morphism of
Lie algebroids} is a smooth vector bundle map $\varphi: A^1\rightarrow A^2$ over
a base map $f:M^1\rightarrow M^2$, satisfying the natural anchor compatibility
condition, $\rho_2\comp\phi=\dd f\comp\rho_1$,
 and a bracket compatibility condition. The bracket compatibility can be
expressed in a few equivalent ways, using decompositions of sections or
connections: pick a Koszul connection $\nabla$ for $A^2\rightarrow M^2$,
denote with $f^{\pmb{+}}\nabla$ the induced connection on the pullback bundle
$f^{\pmb{+}} A^2\rightarrow M^1$ and with $\varphi^!$ the bundle map
$A^1\rightarrow f^{\pmb{+}} A^2$ induced by $\varphi$, the condition is
$$
\varphi^![a,b]_1=f^{\pmb{+}}\nabla_{\rho_1(a)}\varphi^!b -
f^{\pmb{+}}\nabla_{\rho_1(b)}\varphi^!a -f^{\pmb{+}}\tau^\nabla(\varphi^!a,\varphi^!b)
\quad,\quad a,b\in\Gamma(A^1)\qquad,$$
where $f^{\pmb{+}}\tau^\nabla$ is the pullback of the torsion tensor of
$\nabla$. With the above definition one can show that there is a category of Lie
algebroids, with direct products and pullbacks (under natural
transversality conditions).
More invariantly, a Lie algebroid structure on a
vector bundle $A\rightarrow M$ is
equivalent to a diferrential on the graded algebra
$\sf{C}^\bullet(A):=(\Gamma(\wedge^{\bullet}A^*),\wedge)$ \cite{Vain}; a
vector bundle map
$\phi: A^1\rightarrow A^2$  is then a morphism of Lie
algebroids iff \cite{K} the induced map
$\sf{C}^\bullet(A^2)\rightarrow \sf{C}^\bullet(A^1)$ is a chain map.\\
A \textsf{Lie subalgebroid} $B$ of $A$ is a vector subbundle such that the
inclusion $B\hookrightarrow A$ is a morphism of Lie algebroids.\\
The Lie algebroid of a Lie groupoid $\poidd{\calG}{M}$ is the vector
bundle $T^{\sor}_M\calG$ (the restriction of the kernel of the tangent source map to
the base manifold) endowed with a bracket induced from that of right invariant
sections of $T^{\sor}\calG$. Similarly restricting the tangent map of a morphism of
Lie groupoids $\calG^1\rightarrow\calG^2$ to
$T^{\sor^1}_{M^1}\calG^1\rightarrow T^{\sor^2}_{M^2}\calG^2$, one obtains a morphism
of Lie algebroids. That is, there exists a \tsf{Lie functor} from the category of
 Lie
groupoids to that of Lie algebroids.\\
Moreover the Lie functor preserves direct
products and pullbacks; note, however, that fibered products of
Lie groupoids do not always exist.
\begin{proposition}\label{Fibbia}\footnote{Stronger sufficient conditions for the existence of fibred
products of Lie groupoids appeared in \cite{Mb}.}
Consider morphisms of Lie groupoids $\varphi_{1,2}:\calG^{1,2}\rightarrow\calH$
over
$f_{1,2}:M^{1,2}\rightarrow N$, such that $\varphi_1\pitchfork\varphi_2$ and
$f_1\pitchfork f_2$. Then, the manifold fibered product
$\calG^1\fib{\varphi_1}{\varphi_2}\calG^2$ exists and the natural
$($smooth$\,)$ groupoid
structure over $M^1\fib{f_1}{f_2}M^2$ induced from the direct product
$\calG^1\times\calG^2$
is that of a \emph{Lie} groupoid iff
the \emph{\tsf{source transversality condition}}
\begin{equation}\label{@@@}
\dd\varphi_1T^{\sor_1}_{g_1}\calG^1 +
\dd\varphi_2T^{\sor_2}_{g_2}\calG^2=T^{\sor}_h\calH\quad,
\quad \varphi_1(g_1)=h=\varphi_2(g_2)
\end{equation}
is satisfied for all $(g_1,g_2)\in\calG^1\fib{\varphi_1}{\varphi_2}\calG^2$.
In this case, $\calG^1\fib{\varphi_1}{\varphi_2}\calG^2$ is a fibered product in
the category of Lie groupoids.
\end{proposition}
\begin{proof} We have to prove that the source map of
$\calG\equiv\calG^1\fib{\varphi_1}{\varphi_2}\calG^2\to
M\equiv M^1\fib{f_1}{f_2}M^2$ is submersive
iff (\ref{@@@}) holds; universality is then manifest.
Let $q_{1,2}=\sor_{1,2}(g_{1,2})$, $q=\sor(h)$,
for any $(g_1,g_2)\in\calG$, $h=\fii_{1,2}(g_{1,2})$, applying the snake lemma to the exact
commuting diagram
$$
\xy
*+{}="0", <0.5cm,0cm>
*+{0}="Z1", <3cm,0cm>
*+{T^{\sor_1}_{g_1}\calG^1\oplus T^{\sor_2}_{g_2}\calG^2}="S1", <6.5cm,0cm>
*+{T_{g_1}\calG^1\oplus T_{g_2}\calG^2}="G1", <11cm,0cm>
*+{T_{q_1}M^1\oplus T_{q_2}M^2}="T1", <13.5cm,0cm>
*+{0}="ZZ1", <0.5cm,-1.5cm>
*+{0}="Z2", <3cm,-1.5cm>
*+{T^{\sor}_h\calH}="S2", <6.5cm,-1.5cm>
*+{T_h\calH}="G2", <11cm,-1.5cm>
*+{T_qN}="T2",<13.5cm,-1.5cm>
*+{0}="ZZ2"
\ar   @{->} "S1";"S2"^{\Phi_{12}}
\ar   @{->} "G1";"G2"^{\dd\fii_1-\dd\fii_2}
\ar   @{->} "T1";"T2"^{\dd f_1-\dd f_2}

\ar   @{->} "Z1";"S1"
\ar   @{->} "S1";"G1"
\ar   @{->} "T1";"ZZ1"
\ar   @{->} "G1";"T1"^{{}^{{}^{\hbox{\small{$\dd\sor_1\times\dd\sor_2$}}}}}

\ar   @{->} "Z2";"S2"
\ar   @{->} "S2";"G2"
\ar   @{->} "G2";"T2"_{\dd\sor}
\ar   @{->} "T2";"ZZ2"
\endxy
$$
yields a connecting arrow $\pa$ and a long exact sequence
$$\quad
\xy
*+{0}="0", <2cm,0cm>
*+{\mathrm{ker}\,\Phi_{12}}="1", <5cm,0cm>
*+{T_{(g_1,g_2)}\calG}="2", <8cm,0cm>
*+{T_{(q_1,q_2)}M}="3", <11cm,0cm>
*+{\mathrm{coker}\,\Phi_{12}}="4", <13cm,0cm>
*+{0}="5",
\ar   @{->} "0";"1"
\ar   @{->} "1";"2"
\ar   @{->} "2";"3"^{\dd\sor_{12}}
\ar   @{->} "3";"4"^{\pa}
\ar   @{->} "4";"5"
\endxy\quad,
$$
where $\Phi_{12}:=\!\left.(\dd\fii_1-\dd\fii_2)\right|_{
T^{\sor_1}_{g_1}\calG^1\oplus T^{\sor_2}_{g_2}\calG^2}$ and
$\sor_{12}=\!\left.(\sor_1\times\sor_2)\right|_\calG$. The result follows.
\end{proof}
Generalizing slightly the constructions of \cite{HM} allows to introduce
abstract fibered products of Lie algebroids\footnote{The existence of fibered product of Lie
algebroids under the natural transversality conditions was stated without proof in \cite{HM}.}.
Recall from \cite{HM} that the \emph{pullback
algebroid}, $f\daga A\rightarrow N$, along a smooth map $f :N\rightarrow M$
of a Lie algebroid $A\rightarrow M$, with anchor $\rho$,
is defined whenever $\dd f-\rho^!:TN\oplus f^{\pmb{+}} A\rightarrow TM$ has
constant rank (in particular, when $f$ is submersive); the total space of $f\daga A$
is then $\ker\, (\dd f-\rho^!)$, the anchor, denoted with $\rho_{{\pmb{+}}{\pmb{+}}}$,
being the first projection (see {\footnotesize [HM, Mb]} for a description of the bracket). Pullback algebroids satisfy the following
important universal property.
\begin{proposition}[\cite{HM}] Let $A^{(')}\rightarrow M^{(')}$ be Lie algebroids and $f:N\rightarrow M$ a smooth map, such
that the pullback $f\daga A$ exists.
For any morphism of Lie algebroids $\phi: A'\rightarrow A$ over
$g: M'\rightarrow M$, such that there is a smooth factorization $g=f\comp h$,
for some $h:M'\rightarrow N$, there exists a unique morphism $\psi:
A'\rightarrow f\daga A$, such that $\phi=f\daga\comp\psi$, for the natural
morphism $($a.k.a. the \emph{\tsf{inductor}}$)$
$f\daga :f\daga A\rightarrow A$.
\end{proposition}
Next, consider Lie algebroids $A^1$, $A^2$ and $B$ \emph{over the same base}
$M$; given morphisms of Lie algebroids $\phi_{1,2}: A^{1,2}\rightarrow B$ over the identity,
such that the fibered product
$A^1\fib{\phi_1}{\phi_2} A^2$ is a vector bundle, it is possible to introduce the
\tsf{fibered product Lie algebroid} (over $B$ in this case), for the vector bundle structure over $M$. The anchor is
$\rho(a_1\oplus a_2)=\rho_B\comp\phi_1(a_1)
=\rho_B\comp\phi_2(a_2)$, for any
$a_1\oplus a_2\in A_1\fib{\phi_1}{\phi_2} A_2$;
the bracket is defined componentwise:
$$\quad\!\!
\brak{a_1\oplus a_2}{b_1\oplus b_2}=\brak{a_1}{b_1}\oplus
\brak{ a_2}{b_2}\quad,\quad a_1\oplus a_2,b_1\oplus b_2\in
 \Gamma(M,A_1\fib{\phi_1}{\phi_2} A_2)
\quad.$$
Last construction is a straightforward generalization of the \emph{product
of Lie algebroids over the same base} in \cite{HM}, which is recovered
replacing $B$ with $TM$ and $\phi_{1,2}$ with $\rho_{1,2}$.\\
Given
Lie algebroids $A^{1,2}\rightarrow M^{1,2}$, denote with $M^{12}$ the direct product
$M^1\times M^2$ and with $\rm{pr}_{1,2}$ the projections onto $M^{1,2}$. Since
the pullback algebroids $\rm{pr}_{1,2}\daga A^{1,2}$ always exist and the
fibered product of manifolds
$\rm{pr}_1\daga A^1\fib{\rho^{1}_{{\pmb{+}}{\pmb{+}}}}{\rho^{2}_{{\pmb{+}}{\pmb{+}}}}
\rm{pr}_2\daga A^2$ is to be identified with the vector bundle $A^1\times
A^2\rightarrow M^1\times M^2$,
there is always a fibered product Lie algebroid over $TM^{12}$, the
\tsf{direct product Lie algebroid} (denoted
simply as $A^1\times A^2$) of $A^1$ and $A^2$; it is straightforward to check that it is indeed a
direct product in the category of Lie algebroids.
\begin{proposition}\label{fibbia}
Consider morphisms of Lie algebroids $\phi_{1,2}: A^{1,2}\rightarrow B$ over
$f_{1,2}:M^{1,2}\rightarrow N$, such that $\phi_1\,{\pmb{\pitchfork}}\,\phi_2$.
Then, the vector bundle fibered product $A^1\fib{\phi_1}{\phi_2}A^2\rightarrow
M^1\fib{f_1}{f_2}M^2$ carries a unique Lie algebroid structure making it
a Lie subalgebroid of the direct product $A^1\times A^2$, thus  a fibered product in the
category of Lie algebroids.
\end{proposition}
\begin{proof}
Denote with $M$ the fibered product $M^1\fib{f_1}{f_2}M^2$ and
with $p_{1,2}:M\rightarrow M^{1,2}$
the restrictions of the projections on the first and second component.
Transversality for $f_1$ and
$f_2$ implies
$\delta:=\left.(f_1\times f_2)\right|_M=p_{1,2}\comp f_{1,2}$ being submersive
to the diagonal $\Delta_N$; then, upon identifying $N$ with
$\Delta_N$, $B$ can be pulled back to a Lie algebroid over $M$. The pullback
Lie algebroids $p_{1,2}\daga A^{1,2}\rightarrow M$ also exist and there is a
fibered product Lie algebroid $p_1\daga A^1\fib{\psi_1}{\psi_2}p_2\daga A^2$ over
$\delta\daga B$,
where the morphisms $\psi_{1,2}$ are the unique obtained factorizing the
compositions of $\phi_{1,2}$ with the inductors
$p_{1,2}\daga A^{1,2}\rightarrow A^{1,2}$ along the identity of $M$:
$$\qquad
\xy
*+{}="0", <-2.5cm,-0.7cm>
*+{M}="-3",<-2.5cm,0.7cm>
*+{p_{1,2}\daga A^{1,2}}="-1",     <-0.7cm,0.7cm>
*+{A^{1,2}}="1",    <0.7cm,0.7cm>
*+{B}="2", <-0.7cm,-0.7cm>
*+{M^{1,2}}="3", <0.7cm,-0.7cm>
*+{N}="4",   <2cm,-2.1cm>
*+{M}="4'",  <2cm,-0.7cm>
*+{\delta\daga B}="2'"
\ar       @{->} "-1";"2'"^{\psi_{1,2} }
\ar       @{->} "-1";"1"_{}
\ar       @{->} "-1";"-3"_{ }
\ar       @{->} "-3";"3"_{ }
\ar       @{->} "-1";"1"_{ }
\ar       @{->} "1";"3"_{ }
\ar       @{->} "1";"2"^{\phi_{1,2}}
\ar       @{->} "3";"4"^{}
\ar       @{->} "2";"4"_{}
\ar       @{->} "2'";"4'"^{ }
\ar       @{->} "4'";"4"^{}
\ar       @{=} "-3";"4'"^{}
\ar       @{->} "2'";"2"^{}
\endxy\qquad.
$$
Note that $A^1\fib{\phi_1}{\phi_2}A^2$
and $p_1\daga A^1\fib{\psi_1}{\psi_2}p_2\daga A^2$ coincide and are smooth
manifolds due to transversality of $\phi_{1,2}$, thus
$A^1\fib{\phi_1}{\phi_2}A^2$ inherits a Lie algebroid structure from
$p_1\daga A^1\fib{\psi_1}{\psi_2}p_2\daga A^2$.
Next we show that there exists a unique morphism of Lie algebroids
$\chi:A^1\fib{\phi_1}{\phi_2}A^2\rightarrow A^1\times A^2$ filling the
diagram
$$\qquad
\xy
*+{}="0",                <-1.4cm,0.8cm>
*+{A^1\fib{\phi_1}{\phi_2}A^2}="11",    <2.5cm,-0.4cm>
*+{A^1\times A^2}="1", <1cm,0.8cm>
*+{p_2\daga A^2}="22", <4.6cm,-0.4cm>
*+{\mathrm{pr}_2\daga A^2}="2", <-1.4cm,-0.8cm>
*+{p_1\daga A^1}="33", <2.5cm,-2cm>
*+{\mathrm{pr}_1\daga A^1}="3", <1cm,-0.8cm>
*+{TM}="44",  <4.6cm,-2cm>
*+{T M^{12}}="4",
\ar       @{->} "1";"3"_{}
\ar       @{->} "1";"2"^{}
\ar       @{->} "3";"4"_{}
\ar       @{->} "2";"4"^{}
\ar       @{->} "11";"1"_{\chi}
\ar       @{->} "22";"2"^{\chi_2}
\ar       @{->} "11";"22"^{}
\ar       @{->} "33";"3"_{\chi_1}
\ar       @{->} "11";"33"^{}
\ar       @{^(->} "44";"4"^{}
\ar       @{->} "22";"44"^{}
\ar       @{->} "33";"44"^{}
\endxy\qquad.
$$
The maps to $TM$ and $TM^{12}$ are anchors and
$\chi_{1,2}$
$$
\xy
*+{}="0",     <-1.7cm,0.7cm>
*+{p_{1,2}\daga A^{1,2}}="1",    <0.7cm,0.7cm>
*+{A^{1,2}}="2", <-1.7cm,-0.7cm>
*+{M}="3", <0.7cm,-0.7cm>
*+{M^{1,2}}="4",   <2cm,-2.1cm>
*+{M^{12}}="4'",  <2cm,-0.7cm>
*+{\mathrm{pr}_{1,2}\daga A^{1,2}}="2'"
\ar       @{->} "1";"2'"^{\chi_{1,2} }
\ar       @{->} "1";"3"_{ }
\ar       @{->} "1";"2"^{}
\ar       @{->} "3";"4"^{}
\ar       @{->} "2";"4"_{}
\ar       @{->} "2'";"4'"^{ }
\ar       @{->} "4'";"4"^{}
\ar @ <-0.2cm>      @{^(->} "3";"4'"^{}
\ar       @{->} "2'";"2"^{}
\endxy
$$
are the unique
morphisms factorizing the inductors of $p_{1,2}\daga A^{1,2}$
through the inductors of
$\mathrm{pr}_{1,2}\daga A^{1,2}$ along the natural inclusion $M\inc M^{12}$ of the fibered product
base manifold.
From the anchor compatibility condition for $\chi_{1,2}$ and the
commutativity of last diagram, $\chi_{1,2}$  must be the canonical inclusions.
Moreover, post composing $\chi_{1,2}$ with the inductors
$\mathrm{pr}_{1,2}\daga A^{1,2}\rightarrow A^{1,2}$, yields morphisms of Lie algebroids
$p_{1,2}\daga A^{1,2}\rightarrow A^{1,2}$; by universality of $A^1\times A^2$,
$\chi$ exists and is unique. Thus, $\chi$ is the canonical inclusion, due to the commutativity of the
defining diagram. To check universality of $A^1\fib{\phi_1}{\phi_2}A^2$ is now straightforward.
\end{proof}
It is then clear that the Lie algebroid of a fibered product
$\calG^1\fib{\varphi_1}{\varphi_2}\calG^2$ Lie
groupoid coincides, as a vector bundle, with the fibered product
$A^1\fib{\phi_1}{\phi_2}A^2$ for the induced morphisms $\phi_{1,2}$
(the existence
of such being assured from the source transversality condition:
$A^1\fib{\phi_1}{\phi_2}A^2=\ker_{M^1\fib{f_1}{f_2}M^2}(\phi_1 -\phi_2)$,
 where the bundle map $(\phi_1 -\phi_2):\left.A^1\times A^2\right|_{M^1\fib{f_1}{f_2}M^2}
 \rightarrow B$ has maximal rank rank); the induced
Lie algebroid structure is then that of a fibered product, being
$\calG^1\fib{\varphi_1}{\varphi_2}\calG^2\subset \calG^1\times\calG^2$ a
Lie subgroupoid, by uniqueness.
\bgn{example} Consider a Lie algebroid $A\rightarrow M$ and a smooth map
$f:N\rightarrow M$. The requirement $\dd\,f\,\pmb{\pitchfork}\,\rho$ is precisely
 the
transversality condition for the pullback algebroid $f\daga A$ to exist;
in this
case the fibred product Lie algebroid $TN\fib{\dd\,f}{\rho}A\to\gr{\,f}$
also exists and coincides with $f\daga A$ up to the identification
$\gr{\,f}\simeq N$.
\end{example}
\section{A functorial approach to the integration of $\lagpd$s}
\label{fungo}
Following Pradines \cite{P}, a \tsf{differentiable graph} is a pair of manifolds $(\Gamma, M)$, endowed
with surjective submersions $\alpha,\beta:\Gamma\rightarrow M$. We shall say that a graph is
\tsf{unital}, if there is an injective immersion $\eps: M\rightarrow\Gamma$, for which both
$\alpha$ and $\beta$ are left inverses, resp. \tsf{invertible}, if
there is a diffeomorphism $\iota:\Gamma\rightarrow\Gamma$, such that $\iota^2=\id_\Gamma$,
$\alpha\comp\iota=\beta$ and $\iota\comp\eps=\eps$.
Namely, an invertible unital graph is a ``Lie groupoid without a multiplication''. For any
differentiable graph, each \tsf{nerve}
$$
\Gamma^{(n)}=\Gamma\fib{\alpha}{\beta\circ\mathrm{p_1}}
\Gamma^{(n-1)}\subset\Gamma^{\times
n}\quad,
\quad n>1\quad,
$$
is a smooth submanifold (conventionally, $\Gamma^{(1)}=\Gamma$ and $\Gamma^{(0)}=M$).

Let $(\Ohm,A;\calG, M)$ be an $\lagpd$;  if
$\Omega\rightarrow\calG$ is integrable, so is the Lie subalgebroid
$A\rightarrow M$. Denote with $\poidd{\Xi}{\calG}$ and
$\poidd{\calA}{M}$ the source 1-connected integrations. The top
groupoid structural maps
$\sh$, $\th:\Ohm\rightarrow A$, $\epsh:A\rightarrow\Ohm$ and
$\iotah:\Ohm\rightarrow\Ohm$
are morphisms of Lie algebroids and integrate uniquely to morphisms
of Lie groupoids
$\nstv$, $\nttv:\Xi\rightarrow\calA$, $\etv:\calA\rightarrow\Xi$ and
$\itv:\Xi\rightarrow\Xi$.
The compatibility conditions $\sh\circ\epsh=\id_{A}$, 
$\th\circ\epsh=\id_{A}$, 
$\sh\circ\iotah=\th$,  
$\iotah\circ\iotah=\id_{\Ohm}$ and $\iotah\comp\epsh=\epsh$
are diagrams of Lie algebroid morphisms and integrate to analogous relations for
$\itv$, $\nstv$, $\nttv$ and $\etv$. Then $\itv$ is a diffeomorphism (being the inverse to itself), $\etv$ is
injective (it has left inverses), $\nstv$ and $\nttv$ are surjective (being left
inverses to $\etv$). Actually, $\etv$ is an immersion, $\nstv$ and $\nttv$ are
submersive on an open neighbourhood of $\calG$; this can be seen taking the
tangent diagrams to  $\nstv\comp\etv=\id_{\calG}$ and
$\nttv\comp\etv=\id_{\calG}$. The following fact is obvious for Lie groups, by
equivariance under right translations.
\begin{lemma}\label{lss}
Let $\varphi:\calG\rightarrow\calG'$ be a morphism of Lie groupoids over a
submersive base map $f:M\rightarrow M'$, such that $\varphi$ is submersive on
an open neighbourhood of $M$; then $\varphi$ is submersive on all of $\calG$.
\end{lemma}
\begin{proof}
For any $g\in\calG$, define
$\sigma_g:=T_g\calG\cap\dd\fii^{-1}T^{\sor'}_{\fii(g)}\calG'$, fix an arbitrary
complement $\kappa_g$ and note that $T_g^{\sor}\calG$,
$\ker_g\dd\fii\subset\sigma_g$;
thus $\kappa_g\simeq\dd\fii\,\kappa_g$, $
T^{\sor'}_{\fii(g)}\calG'\cap\dd\fii\,\kappa_g=\{0\}$ and
it is sufficient to show:
(1) $\dd\fii: T^{\sor}_{g}\calG\rightarrow T^{\sor'}_g\calG$ is surjective,
(2) $\sf{dim}\,\kappa_g=\sf{dim}\,M'$.
To see (1), for any $\delta'\in T^{\sor'}_{\fii(g)}\calG'$ pick a
$\delta_o\in T_{\eps(\tar(g))}\calG\cap\dd\fii^{-1}(\dd r_{\fii(g)}^{-1}\delta')$, then
$\delta:=\dd r_g(\delta_o-\dd\eps\,\dd s\, \delta_o)\in T^{\sor}_{g}\calG$ maps on
$\delta'$, by equivariance under right translation of the restriction of
$\dd\fii$ to the subbundle of vectors tangent to the source fibres.  To see (2),
consider that, picking a local bisection $\Sigma_g$ of $\calG$ through $g$,
yields a splitting $\sigma_g=T^{\sor}_g\calG\oplus\dd\Sigma_g\ker_{\sor(g)}\dd f$ and
count dimensions. Alternatively, it is sufficient to check (1) and apply the
5-lemma to the diagram for $\phi$, $f$ and the short exact sequences associated
to the tangent source maps.
\end{proof}

It follows that  $\nstv$ and $\nttv$ are surjective submersions and
$(\Xi,\calA)$ is an invertible unital graph in the category of Lie groupoids.
Moreover, the top vertical nerves of the graph
$$
\qquad
\Xi^{(n)}_V=\Xi\fib{\stv}{\ttv\circ\mathrm{pr_1}}
\Xi^{(n-1)}_V\subset\Xi^{\times
n}\qquad,
\qquad n>1\qquad,
$$
($\Xi^{(1)}_V=\Xi$, $\Xi^{(0)}_V=\calA$) can be inductively endowed
with Lie groupoid structures (transversality conditions are met) over
the side vertical nerves
$\calG^{(n)}$, the groupoid being induced by $\poidd{\Xi}{\calG}$.

\medskip In general, the following holds.
\begin{theorem}
Let $A^{1,2}$ and $B$ be integrable Lie algebroids; consider
morphisms of Lie algebroids $\phi_{1,2}:A^{1,2}\rightarrow B$, such
that transversality conditions of proposition $(\ref{fibbia})$ hold and denote with
$\varphi_{1,2}:\calG^{1,2}\rightarrow \calH$ any integrations. Then, the fibred
product Lie groupoid $\calG^1\fib{\varphi_1}{\varphi_2}\calG^2$  exists.
\end{theorem}
\begin{proof}
With the same notations as in \S\ref{3}, consider that transversality for
$\phi_{1,2}$ is equivalent to the source transversality condition (\ref{@@@})
for $\fii_{1,2}$ along $M$; to see that the condition holds off the base
manifold, use a right translation argument as in lemma (\ref{lss}). It is now
sufficient to prove transversality for $\fii_{1,2}$, i.e to showing that
$(\dd\fii_1 +\dd\fii_2):T_{g_1}\calG^1\oplus T_{g_2}\calG^2\rightarrow T_h\calH$
is surjective for all $(g_1,g_2)\in\calG^1\fib{\varphi_1}{\varphi_2}\calG^2$;
clearly the property holds sourcewise, thus, set
$$\qquad\sigma_{(g_1,g_2)}:=T_{g_1}\calG^1\oplus T_{g_2}\calG^2\cap
(\dd\fii_1 +\dd\fii_2)^{-1}T^{\sor}_h\calH\qquad.$$
As in the proof of lemma (\ref{lss}),
$\ker_{(g_1,g_2)}(\dd\fii_1 +\dd\fii_2)$,
$T^{\sor^1}_{g_1}\calG^1\oplus T^{\sor^2}_{g_2}\calG^2\subset\sigma_{(g_1,g_2)}$ and
for any choice of bisections $\Sigma^{1,2}_{g_{1,2}}$ of $\calG^{1,2}$,
there is a splitting
$$\qquad
\sigma_{(g_1,g_2)}=(T^{\sor^1}_{g_1}\calG^1\oplus T^{\sor^2}_{g_2}\calG^2)
\oplus
(\dd\Sigma^1_{g_1}\oplus\dd\Sigma^2_{g_2})\ker_{(g_1,g_2)}(\dd f_1 +\dd
f_2)\qquad;
$$
to conclude the proof, count dimensions, using transversality of the base maps $f_{1,2}$, to show that
$(\dd\fii_1 +\dd\fii_2)$ maps any complement of $\sigma_{(g_1,g_2)}$
 to a complement of
$T^{\sor}_h\calH$. As in lemma (\ref{lss}) one can replace the last part of this
proof applying the 5-lemma to the suitable exact commuting diagram.
\end{proof}
Actually, the graph $\Xi$ carries more structure. Let $\mathbb{S}_\Xi$ denote the double source
map $(\nsth,\nstv):\Xi\rightarrow\calG\fib{\ssh}{\ssv}\calA$; submersivity for $\mathbb{S}_\Xi$ is
equivalent to fibrewise surjectivity of the restriction
$\left.\dd\nstv\right|_{T^{\sth}\Xi}:T^{\sth}\Xi\rightarrow
T^{\ssh}\calA$, which can be seen to hold by the usual right translation argument and
fibrewise surjectivity of $\sh:\Ohm\rightarrow A$, that is, surjectivity of the double source
map $\$$.
Then $(\Xi,\calG,\calA, M)$ fails to be a double Lie groupoid only for lacking a top
vertical multiplication. Note however, that the Lie algebroids of the top vertical nerves
are the nerves
$$
\Hat{\Ohm}^{(n)}=\Ohm\fib{\sh}{\th\circ\mathrm{pr_1}}
\Hat{\Ohm}^{(n-1)}\subset\Ohm^{\times
n}\qquad, \qquad n>1
$$
of $\poidd{\Ohm}{A}$ and the associativity diagram
$$
\xy
*+{}="0",    <-1cm,1cm>
*+{\Hat{\Ohm}^{(3)}_V}="1", <1cm,1cm>
*+{\Hat{\Ohm}^{(2)}_V}="2", <-1cm,-0.7cm>
*+{\Hat{\Ohm}^{(2)}_V}="3", <1cm,-0.7cm>
*+{\Ohm}="4",
\ar       @{->} "1";"3"_{\muh\times\id_\Ohm}
\ar       @{->} "1";"2"^{\id_\Ohm\times\muh}
\ar       @{->} "3";"4"_{\muh}
\ar       @{->} "2";"4"^{\muh}
\endxy
$$
commutes in the category of Lie algebroids; thus, whenever the second
top vertical nerve of $\Xi$ is source 1-connected, $\muh$
integrates to a morphism
$\mtv:\Xi^{(2)}_V\rightarrow\Xi$ and the diagram
$$
\xy
*+{}="0",    <-1cm,1cm>
*+{\Xi^{(3)}_V}="1", <1cm,1cm>
*+{\Xi^{(2)}_V}="2", <-1cm,-0.7cm>
*+{\Xi^{(2)}_V}="3", <1cm,-0.7cm>
*+{\Xi}="4",
\ar       @{->} "1";"3"_{\mtv\times\id_\Xi}
\ar       @{->} "1";"2"^{\id_\Xi\times\mtv}
\ar       @{->} "3";"4"_{\mtv}
\ar       @{->} "2";"4"^{\mtv}
\endxy
$$
commutes in the category of Lie groupoids, provided the third top vertical nerve
is source connected. To see this, consider that there is a similar
\emph{commuting} diagram, where $\Xi^{(3)}_V$ is replaced by its source
1-connected cover; commutativity of the associativity diagram follows by
diagram chasing due to surjectivity of the covering morphism.
Compatibility of $\mtv$ with the graph over
$(\Xi,\calA)$ can be easily shown by integrating
 the suitable compatibility diagrams for $\muh$ to be a groupoid multiplication.

Under source connectivity assumptions on the second and third top vertical
nerves, also integrability results for morphisms and sub-objects follow in the
same fashion.

\bigskip

Even though no counter-examples are known to us, it seems unlikely
for $\Xi^{(2,3)}_V$ to be source (1-)connected in general, for the vertically
source 1-connected graph $\Xi$ of an $\lagpd$. A sufficient condition for that to hold
is:\medskip\\
\emph{for each $g\in\calG$, the restriction
$\nstv:\nsth^{-1}(g)\rightarrow\sor_{h}^{-1}(\sor_v(g))$ has the 1-homotopy lifting
property}.
\medskip\\ Last requirement can be reduced to a lifting condition
with respect to post composition with $\sh$ for $A$-paths,
resp. $A$-homotopies, to
$\Ohm$-paths, resp. $\Ohm$-homotopies (regarding $\Ohm$ as a Lie algebroid over $\calG$),
involving only the $\lagpd$-data.\\
Recall from \cite{CrFr}, that, for any Lie algebroid $(A\rightarrow M, [\:,\:], \rho)$, an $A$-path is a
$\calC^1$ map $\alpha:I\rightarrow A$, over a base path
$X:I\rightarrow M$, such that $\dd X=(\rho\comp X)\, \alpha$, namely a morphism of Lie algebroids
$TI\rightarrow A$. An $A$-homotopy $h$ from an $A$-path $\alpha_+$ to an $A$-path $\alpha_-$ is a
$\calC^1$ morphism of Lie algebroids $TI^{\times 2}\rightarrow A$, satisfying suitable boundary
conditions. Regarding $h$ as a 1-form taking values in the pullback $Y^{\pmb{+}} A$, for the base map
$Y$, the anchor compatibility condition is $\dd Y=(\rho\comp Y)\, h$, while the bracket compatibility takes
the form of a Maurer-Cartan equation
$$\qquad
\mathrm{D}\,h+\frac{1}{2}[h
\,\overset{\wedge},\,
h]=0\qquad,
$$
where $\mathrm{D}$ is the covariant derivative of the pullback of an arbitrary connection $\nabla$
for $A$ and $
[\theta_+\!
\,\overset{\wedge},\,
\theta_-]$ is the contraction $-\iota_{\theta_+\wedge\theta_-}Y^{\pmb{+}}\tau^\nabla$ with the pullback
of the torsion tensor, $\theta_\pm\in\Ohm^1(Y^{\pmb{+}} A)$. The boundary conditions are
$\iota^*_{\pa^\pm_{\mathrm{H}}}h=\alpha_\pm$ and
$\iota^*_{\pa^\pm_{\mathrm{V}}}h=0$, where
$$
\qquad
\pa^-_{\hbox{\tiny{$\mathrm{H(,V)}$}}}=\{(x,y)|x(y)=0\}
\quad\hbox{ and }\quad
\pa^+_{\hbox{\tiny{$\mathrm{H(,V)}$}}}=\{(x,y)|x(y)=1\}\qquad.
$$
Assume $A$ is integrable and $\calG$ the source 1-connected integration; then, the space of class
$\calC^2$ $\calG$-paths and of $A$-paths
are homeomorphic for the natural Banach topologies, the homeomorphism
$\calG\hbox{-}paths\rightarrow A\hbox{-}paths$ being given by the ``right
derivative'' $\delta\,\Xi(t)=\dd r_{\Xi(t)}^{-1}\dot{\Xi}(t)$, $t\in I$. Under such
correspondence homotopy of $\calG$-paths within the source fibers and relative to the endpoints, viz.
$\calG$-homotopy,
translates precisely into $A$-homotopy.

\medskip
We are ready to state our conditions.
\begin{proposition}\label{lift}
With the same notations as above: all the nerves $\Xi^{(n)}_V$, $n\in\NN$,
of $\Xi$ are:\smallskip\\
$i)$ Source connected, if $\sh$ has the
\emph{\tsf{0-$\mathcal{LA}$-homotopy lifting property}},
i.e.
for any $\Ohm$-path $\xi$ and $A$-path $\alpha$,
which is $A$-homotopic to $\sh\comp\xi$, there exists a $\Ohm$-path $\xi'$,
which is $\Ohm$-homotopic to $\xi$ and satisfies
$\sh\comp\xi'=\alpha$, resp.\smallskip\\
$ii)$ Source 1-connected, if $\sh$ has the
\emph{\tsf{1-$\mathcal{LA}$-homotopy lifting property}},
i.e.
for any $\Ohm$-path $\xi$, which is $\Ohm$-homotopic to the constant
$\Ohm$-path $\xi_o\equiv 0^\Ohm_{\mathrm{Pr}(\xi(0))}$,
and
$A$-homotopy $h$ from  $\alpha:=\sh\comp\xi$ to the constant $A$-path  $\alpha_o\equiv 0^A_{\mathrm{pr}(\alpha(0))}$,
there exists a $\Ohm$-homotopy $\hat{h}$ from $\xi$ to $\xi_o$,
such that $\sh\comp\hat{h}=h$.
\end{proposition}
\begin{proof} We shall prove the statement for $n=2$, the argument extends easily by
induction. ($i$) For any $x_\pm\in\Xi^{(2)}_V$ we have to find $\Xi$-paths $\gamma_\pm$,
such that $\gamma_\pm(1)=x_\pm$ and $\nstv\comp\gamma_+=\nttv\comp\gamma_-$. For any
$\Ohm$-paths $\xi^\pm\in x_\pm$, $\alpha_+:=\sh\comp\xi^+$ is $A$-homotopic to
$\alpha_-:=\th\comp\xi^-$ and we can assume $\alpha_-=\sh\comp\xi^+$. Define smooth families of
$\Ohm$-paths $\xi^\pm_t(u):=t\cdot\xi^\pm(t u)$, $u$, $t\in I$; the corresponding (unique) families
$\{\omega^\pm_t\}_{t\in I}$ of $\Xi$-paths are
contained within the source fibers of $x_\pm$:
$$\qquad
\nsth(\omega^\pm_t(u))=\epsh(\Pr(\xi^\pm_t(0)))=\epsh(\Pr(\xi^\pm(0)))=
\nsth(x_\pm)\qquad.
$$
Thus, projecting to the quotient ${\cal W }(\Ohm)$, yields $\Xi$-paths
$\gamma_\pm(t)=[\xi^\pm_t]$,
with $\gamma_\pm(1)=[\xi^\pm]=x_\pm$ and
$$\qquad\nstv(\gamma^+(t))=[\nstv\comp\omega^+_t]=[\sh\comp\xi^+_t]=
[\th\comp\xi^-_t]=[\nttv\comp\omega^-_t]=\nttv(\gamma^-(t))\qquad.$$
($ii$) For any $\Xi^{(2)}_V$-loop $(\lambda^+,\lambda^-)$,
we have to find $\Xi$-homotopies $H_\pm$ from $\lambda^\pm$ to the constant $\Xi$-loops
$\lambda^\pm_o\equiv\lambda^\pm(0)$, such that $\nstv\comp H^+=\nttv\comp H^-$. We may as well regard
$\lambda^\pm$ as $\Xi$-loops, satisfying $\nstv\comp \lambda^+=\alpha=\nttv\comp\lambda^-$,
for some $\cal A$-loop $\alpha$; thus,
the (unique) corresponding $\Ohm$-paths $\xi^\pm$ satisfy $\sh\comp\xi^+=\hat{\alpha}=\th\comp\xi^-$,
for the $A$-loop $\hat{\alpha}$ corresponding to $\alpha$. Any $\Ohm$-homotopy $h^-$ from
$\xi^-$ to the constant $\Ohm$-path $\xi^-_o$ projects down to an $A$-homotopy $\th\comp h^-$
from $\hat{\alpha}$ to the constant $A$-path $\hat{\alpha}_o$ and can be lifted along $\sh$
to an $\Ohm$-homotopy
$h^+$ from $\xi^+$ to $\xi^+_o$, since $\Xi$ is 1-connected and $\xi^+$ the $\Ohm$-path of
a $\Xi$-loop, therefore $\Ohm$-homotopic to the constant $\Ohm$-path $\xi^+_o$.
The $\Xi$-homotopies $H_\pm$ corresponding to $h^\pm$ are defined by
$H^\pm(\eps,t)=\tilde{h}_\pm(\eps,t,0,0)$, $\eps$, $t\in I$, where
$\Tilde{h}^\pm:I^{\times 2}\times I^{\times
2}\rightarrow\Xi$ are the morphisms of Lie groupoids integrating
$h^\pm:TI^{\times 2}\rightarrow\Ohm$. Therefore the condition $\nstv\comp H^+=\nttv\comp H^-$ follows
integrating $\sh\comp h^+=\th\comp h^-$ and $H_\pm$ are the desired homotopies.
%
%
%
%
%
%
%
%
%
%
%
%
%
%
%
%
%
%
\end{proof}
Therefore, we have:
\bgn{theorem}
Any $\lagpd$ $\sf{\Ohm}$ with integrable top Lie algebroid,
 whose associated integrating graph
has source 1-connected top vertical nerve and source connected third
top vertical nerve, is integrable to a horizontally source 1-connected double Lie
groupoid. The connectivity assumptions hold when the
$\mathcal{LA}$-homotopy lifting conditions of proposition $(\ref{lift})$ on
$\sf{\Ohm}$ are met.
\end{theorem}
\begin{example} Consider the tangent prolongation $\lagpd$ of example \ref{exa}. For any
manifold $N$, $TN$-paths, resp. $TN$-homotopies, may be identified with ordinary paths in
$N$, resp. homotopies in $N$ relative to the endpoints. If the source map of $\calG$ has the
1-homotopy lifting property\footnote{Necessary and sufficient conditions for a surjective
submersion to be an $n$-fibration, $n\in\NN$, can be found in \cite{Mg}.}, in particular when
$\poidd{\calG}{M}$ is a proper groupoid, our conditions are satisfied and the
fundamental groupoid
$\poidd{\Pi(\calG)}{\calG}$ carries a further Lie groupoid structure over $\Pi(M)$,
making
$$
\xy
*+{}="0",    <-0.7cm,0.7cm>
*+{\Pi(\calG)}="1", <0.7cm,0.7cm>
*+{\calG}="2", <-0.7cm,-0.7cm>
*+{\Pi(M)}="3", <0.7cm,-0.7cm>
*+{M}="4",
\ar  @ <0.07cm>   @{->} "1";"2"^{}
\ar  @ <-0.07cm>  @{->} "1";"2"_{}
\ar  @ <-0.07cm>  @{->} "1";"3"_{}
\ar  @ <0.07cm>   @{->} "1";"3"^{}
\ar  @ <0.07cm>   @{->} "2";"4"^{}
\ar  @ <-0.07cm>  @{->} "2";"4"_{}
\ar  @ <-0.07cm>  @{->} "3";"4"_{}
\ar  @ <0.07cm>   @{->} "3";"4"^{}
\endxy
$$
a double Lie groupoid; the top vertical multiplication is induced by
multiplication of
pointwise composable paths in $\calG$. The double Lie groupoid of example
(\ref{exa}) is
clearly an integration of the same $\lagpd$ (even when $\calG$ is not proper),
which is, in general, not horizontally
source 1-connected; if $\calG$ itself is source 1-connected, since $\sor$ is a
submersion with 1-connected fibers, it is possible to apply proposition (7.1) of
\cite{Ginz} to conclude that $\pi_{0,1}(\calG)=\pi_{0,1}(M)$. Thus,
when $\pi_{0,1}(M)=0$, our approach produces the pair double groupoid of $\cal{G}$.
\end{example}
Consider any double Lie groupoid;
since there is no natural way to to replace the horizontal
Lie groupoids with the 1-connected
covers of their source connected components without affecting the top vertical groupoid,
in general, one should not expect to be able to integrate an $\lagpd$ to a horizontally source
1-connected double Lie groupoid.
 However, examples in which our approach applies
should arise from proper actions of Lie groupoids and Lie algebroids;
we shall study elsewhere less restrictive conditions
to integrate $\muh$ to a groupoid multiplication for $(\Xi,\calA)$ and discuss more
interesting examples.
\section{The case of an integrable Poisson groupoid}\label{pioppo}
Consider a Poisson groupoid $\pgpd$, with an integrable Poisson structure and
the associated $\lagpd$ (as in \S\ref{lapo}); denote
with $\calS$ the source 1-connected symplectic groupoid integrating $(\calP,\Pi)$ and
with $\poidd{(\ol{\calP},\ol{\Pi})}{M}$ the source 1-connected weak dual Poisson groupoid.
The graph of $\calS$ over $\ol{\calP}$ is obtained integrating the structural
maps of the graph underlying the cotangent prolongation groupoid
$\poidd{\tsp}{A^*}$. The next result follows specializing the results of last
section.
\bgn{theorem}\label{ddd}
Let $\pgpd$ be a Poisson groupoid with integrable Poisson structure. Then, the
graph
\bgn{equation}\label{dd}
\xy
*+{}="0",    <-0.7cm,0.7cm>
*+{\calS}="1", <0.7cm,0.7cm>
*+{\calP}="2", <-0.7cm,-0.7cm>
*+{\ol{\calP}}="3", <0.7cm,-0.7cm>
*+{M}="4",
\ar  @ <0.07cm>   @{->} "1";"2"^{}
\ar  @ <-0.07cm>  @{->} "1";"2"_{}
\ar  @ <-0.07cm>  @{->} "1";"3"_{}
\ar  @ <0.07cm>   @{->} "1";"3"^{}
\ar  @ <0.07cm>   @{->} "2";"4"^{}
\ar  @ <-0.07cm>  @{->} "2";"4"_{}
\ar  @ <-0.07cm>  @{->} "3";"4"_{}
\ar  @ <0.07cm>   @{->} "3";"4"^{}
\endxy
\end{equation}
integrating the $\lagpd$ associated with $\calP$ is a symplectic double
groupoid, whenever it has source 1-connected second vertical nerve and
source connected third vertical nerve, this is the case when the
$\mathcal{LA}$-homotopy lifting conditions of proposition $(\ref{lift})$
are met.
Moreover, $\calS$ is a double of $\calP$.
\end{theorem}
Note that we make no connectivity assumptions on $\calP$.
\bgn{proof}
Under the connectivity assumptions, the procedure of
\S\ref{fungo} yields a horizontally source 1-connected double
Lie groupoid on (\ref{dd}),
whose top horizontal structure is the symplectic groupoid of $(\calP,\Pi)$.
Recall that the (ordinary) graph of the cotangent multiplication $\muh$ is
(by construction indeed \cite{CDW})
$$
\qquad
\gr{\,\muh}=(\id_{\tsp}\times\id_{\tsp}\times -\id_{\tsp})\, N^*\gr{\,\mu}
\subset\tsp^{\times 3}
\qquad,
$$
where $\mu$ is the
multiplication of $\calP$, therefore it is a Lagrangian submanifold with respect to the
symplectic form
$\ohm_{can}\times\ohm_{can}\times -\ohm_{can}$.
Moreover,
$\gr{\,\mu}$ is coisotropic in
$(\calP^{\times 3},\Pi\times\Pi\times -\Pi)$, thus
$\gr{\,\muh}\rightarrow\gr{\,\mu}$ is canonically a Lie subalgebroid of
$\tsp^{\times 3}\rightarrow \calP^{\times 3}$ for the Koszul bracket. From
integrability of $\Pi$, it follows that $\gr{\,\muh}$ integrates to an
\emph{immersed} Lagrangian \cite{C} subgroupoid $\calL\subset\calS^{\times 3}$
(symplectic form with a minus sign on the third component). Note that the graph
$\gr{\,\varphi}$ of a morphism of Lie groupoids
$\varphi:\calG\rightarrow\calG'$ is always a Lie
subgroupoid of $\calG\times\calG'$ over the graph of its base map,
isomorphic to the domain groupoid. Since the projection to the first and second factor
$\gr{\,\mtv}\rightarrow\calS^{(2)}_V$ is an isomorphism of Lie groupoids,
$\gr{\,\mtv}$ is also source a 1-connected Lie groupoid,
therefore coinciding with the Lie
groupoid $\calL$ of $\gr{\,\muh}$, and
a Lagrangian \emph{Lie} subgroupoid of $\calS^{\times 3}$. Hence $\calS$ is a
symplectic double groupoid.\\
Moreover, the unique Poisson structure $\Pi'$ induced by $\calS$ on $\ol{\calP}$,
makes $\poidd{(\ol{\calP},\Pi')}{M}$ a Poisson groupoid weakly dual to $\pgpd$
and $\Pi'=\ol{\Pi}$, by uniqueness (theorem (\ref{intb})); that is,
$\calS$ is a double for $\pgpd$.
\end{proof}
In the proof of last theorem the fact that the second vertical nerve of the
graph on $\calS$ is 1-connected is not only essential to define the top
horizontal multiplication, but also to show that it is compatible with the
symplectic form. Note that $\poidd{\calS}{\ol{\calP}}$, will not be,
in general, source 1-connected.


\bigskip\bigskip

We conclude this note remarking that Lu and Weinstein's approach
to the construction of the double of a Poisson group does not apply to Poisson groupoids, since it
relies on the peculiar properties of the coadjoint action of a Poisson group on its dual, features
lacking in the theory of Lie groupoids. Moreover the Drinfel'd doubles of Lie bialgebroids are no
longer Lie algebroids. In fact, they can be characterized as double Lie algebroids (Mackenzie's
double), or \cite{ZWX} Courant algebroids.

On the other hand, to any Courant algebroid, one can associate a tridimensional topological field
theory, named
\emph{Courant sigma model} by Strobl \cite{S}. The action functional of this theory was obtained by Ikeda
as a deformation in the BRST antifield formalism
of the abelian Chern-Simons theory coupled with a BF term, but it can also be understood as the
AKSZ action for a symplectic $Q$-manifold of degree 2 \cite{R}.

The phase space of the classical
version of this theory can always be described as constrained hamiltonian system. In the case of the
Courant algebroid of a Lie bialgebroid $(A,A^*)$, with integrable Lie algebroids,
choosing the square as a source
manifold and suitable boundary conditions, the reduced space, when a manifold, carries a symplectic
form and two compatible differentiable (unital, invertible) graphs over the Poisson groupoid of
$(A,A^*)$ and its weak dual. Moreover the geometry of the source manifold allows to define two
partial multiplications before taking the symplectic quotient and yields two compatible Lie groupoid
multiplications on the reduced space, provided some obstruction related to the $\Pi_2$'s of the
leaves of $A$ and $A^*$ vanishes. The chance of integrating simultaneously a Poisson groupoid and its weak dual
by  taking a quotient of the reduced space of the Courant sigma model
is under investigation in a joint work with Alberto Cattaneo.

\noindent
      Luca Stefanini \\
      Institut f\"ur Mathematik, \\
      Universit\"at Z\"urich-Irchel, \\
      Winterthurerstrasse 190, \\
      CH-8057 Z\"urich, \\
      Switzerland \\
      \tsf{lucaste@math.unizh.ch}\\ and\\ \tsf{cucanini@gmail.com}

\label{lastpage}\newpage\markboth{}{} $ $
\end{document}